\documentclass[12pt]{amsart}

\begin{document}

\begin{center}
{\bf  Functions measuring smoothness and the constants in
Jackson--Stechkin theorem }
\end{center}

\vskip .3cm
\begin{center}
A.G. Babenko, Yu.V. Kryakin
\end{center}
\vskip .3cm

\centerline{\em Dedicated to Professor Vitaly Andrienko on the
occasion of his  ${70}^{\text{th}}$ birthday}

\vskip 1cm
\begin{center}
{\bf 1 Introduction }
\end{center}
\vskip .2cm

This paper is devoted to the equivalence of two type direct
theorems in Approximation Theory:

\vskip .2cm

a) for smooth functions  (Favard's estimates ).

\vskip .2cm

b) for arbitrary continuous function (Jackson--Stechkin
estimates).

\vskip .2cm

Specifically, we will show that Jackson--Stechkin inequality with
optimal  respect to the order of smoothness constants follows
from Favard's inequality.

The main tool for this is the function $W_{2k}$, measuring the
smoothness of  integrable periodic function. This characteristic
is more delicate than standard modulus of continuity of the
$2k$-th order. The function  $W_{2k}$ allows us to obtain asymptotically
sharp results for approximation by Favard-type operators. For
example, we obtain the Jackson--Stechkin inequality for periodic
splines with  optimal constants.

Two facts play a key role here.

1. Uniform (on $k$)  boundedness of the operators  $W_{2k}$:

$$  W_{2k}(f,\delta)  \le 3 \| f \|. $$

\vskip .2cm

2. Bernstein--Nikolsky--Stechkin inequality in terms of $W_{2k}$.

\vskip .2cm

This paper is organized as follows. In the Second Section   we
consider the smooth characteristic $W_{2k}$ and prove the
uniformly boundedness of $W_{2k}$. Secton 3 is devoted to
Bernstein--Nikolsky--Stechkin estimate (Theorem 1). Main result of
the paper (Theorem 2) claims that the Favard operators gives the
Jackson--Stechkin theorem with optimal constants. This result is
the consequence of the sharp inequality for the trigonometric
approximation (see \cite{fks}) and will be present in the Forth
Section.

\newpage

\begin{center}
{\bf 2 Functions measuring smoothness }
\end{center}

\vskip .2cm

Let the function  $f$ be continued on the one-dimensional torus $
T=R/ (2 \pi Z)$. The standard smooth characteristic of $f$ is the
modulus of continuity of $r$--th order:

$$ \omega_r (f, \delta):= \sup_{|h| \le \delta} \|\Delta_h^r f \|
= \sup_{|h|\le \delta}  \sup_x |\Delta_h^r f(x) |, $$
where
 $$ \Delta_h^r f(x):= \sum_{j=0}^r (-1)^j \binom{r}j f(x +jh).
$$
We will construct the operators  $W_{2k}$ on the base of the even central difference
$$
\widehat\Delta_h^{2k} f(x):= \sum_{j=-k}^k (-1)^j \binom{2k}{k+j} f(x
+jh).
$$
Define

$$ W_{2k}(f,x,h):= {\binom{2k}k}^{-1} \int_T
\widehat{\Delta}_{t}^{2k} f(x) \phi_h(t) \, dt, \quad  0<h< \pi/k,
$$
where
$$
\phi_h(t) = \begin{cases} \frac 1h (1-\frac{|t|}h), \ &t \in [-h,h], \\
0, & t \notin [-h,h].
\end{cases}
$$
Put
$$
W_{2k}(f, h):=\| W_{2k}(f,\cdot,h) \|.
$$
Introduce two notations. These notations corresponds to sharp maximal function
and maximal function of operator  $W_{2k}$.

$$
W_{2k}^{\sharp} (f,x,\delta):= \sup_{ 0 < h \le \delta} \left| W_{2k}
(f,x,h)\right|.
$$

$$
W_{2k}^{*} (f,\delta):= \sup_{x \in T} W_{2k}^{\sharp}
(f,x,\delta).
$$
One may consider the operators  $W_{2k}$  not only for continuous functions.
The boundedness at all points and integrability  will  suffice.

\
\

\vskip .2cm
{ \bf Properties of the operators  $ W_{2k} , \ W_{2k}^{\sharp}, \  W_{2k}^{*} $}
\vskip .2cm

{\bf 2.1.} The functions $ W_{2k}^{\sharp} (f,x,\delta) \le  \
W_{2k}^{*} (f,\delta)$ are non increasing, as the functions of~$\delta$

\vskip .2cm

{\bf 2.2.}
$$
W_{2k}(f, \delta)  \le W_{2k}^{*} (f, \delta) \le {\binom{2k}k}^{-1} \omega_{2k}(f,\delta).
$$

\vskip .2cm {\bf 2.3.} For the functions that are orthogonal to a
space of trigonometric polynomials of degree  $ \le n-1$ (notation
$g \in T_{n-1}^\perp $) we have \cite{fks}

$$ W_{2k}(g,\delta)  \approx \| g \|, \quad \delta = \alpha \pi/n,
\ \alpha \in (1, n/{k}), $$ or more precisely $$ \| g \| \le
c_\alpha \   W_{2k} (g, \alpha \pi/n ) \le \ c_\alpha  \,
( 1 + \pi^2/8)  \   \| g \|, \eqno (1) $$ where $$ c_\alpha
\le  \sec({\pi}/(2\alpha) ) \le \frac 4\pi \left( 1 - \alpha^{-2}
\right)^{-1}. $$
 \vskip .2cm
\noindent
In the case  $\alpha =1$ we have $$
W_{2k} (g, \pi/n) \approx^{\sqrt{k}} \| g \|,  $$ or, in the
explicit form $$ \| g \| \le  c \ \sqrt{2k} \  W_{2k} (g,
\frac{\pi}{n} ) \le \ c  \ (6/\pi) \ \sqrt{2k}  \    \| g \| ,
\eqno (2) $$ and the first inequality is sharp with respect to
$k$.

\vskip .2cm

The estimates (1), (2)  are the key estimates of  \cite{fks}, devoted to
Jackson--Stechkin inequality with asymptotically sharp constants.

\vskip .2cm
{\bf 2.4.} {\bf Lemma 1.} { \it

$$
W_{2k}^{*} (f, h) \le 3 \| f \|, \quad  h \in (0, \pi/k).
$$
}

\
\
\
\
{\it Proof.} One can rewrite the  function  $W_{2k} (f,x,h)$  in the following form

$$
W_{2k} (f,x,h) = f(x) + (f*\Lambda_{k,h})(x).
$$
Here
$$
\Lambda_{k,h}(x)= 2\sum_{j=1}^k (-1)^{j+1} a_j \phi_{jh}(x), \quad
\phi_{jh}(x):= \frac 1j \phi_h\left(\frac xj \right), \quad a_j:=
\frac{\binom{2k}{k+j}}{\binom{2k}k}.
$$
For estimate of the convolution it is sufficient to put  $h=1$. In this case
$\Lambda_k$ is the even, piecewise-linear function with the vertexes in the points  $(i,b_i), \ i=-k,\dots, k$,
$$
b_{-i}=b_{i}, \quad i=0,\dots, k-1, \quad b_{-k}=b_k=0,
$$
and
$$
b_i= 2 {\binom{2k}k}^{-1} \sum_{j=i+1}^k \binom {2k}{k-j}
(-1)^{j+1} \frac 1j \left(1 - \frac ij \right), \quad i=0,\dots
k-1.
$$
The inequalities
$$
0<b_0< 2 \ln 2, \quad 0>b_1>2\ln2 - \pi^2/6, \quad |b_i|< \frac
1{2i^2}, \quad i=2, \dots k-1,
$$
imply
$$
\int_R |\Lambda_k (t) | \, dt \le 2 \left(\int_0^2 +
\int_2^\infty \right) \le 1.5 + 0.5 = 2.
$$

\qed

\vskip .4cm

\begin{center}
{\bf 3  Bernstein--Nikolsky--Stechkin inequality }
\end{center}
\vskip .2cm

The Bernstein--Nikolsky--Stechkin inequality (see \cite{kbl}, Theorem 3.1.4 in russian edition ) is the generalization
of the classical Bernstein's inequality for trigonometric
polynomials $ \tau \in T_{n}$:

$$
\| D^r \tau  \| \le n^r \| \tau \|,
$$

and reads as

$$
\|D^r \tau \| \le n^r  (2\sin(nh/2))^{-r}  \| \widehat\Delta_h^r \tau \|,
\quad h \in (0, 2\pi/n).
$$

Note, that quantity  $i^r (2\sin(nh/2))^{r}$ is the proper value of the operator
 $\widehat\Delta_h^{r}$ with respect to the eigenfunction $ \exp (int)= c_n(t)+ i s_n(t)$.

\vskip .2cm

{\bf Theorem 1.}{ \it If $\tau \in T_n$,  then

$$
\| D^{2k} \tau \| \le n^{2k}
W_{2k} (c_n,  h)^{-1}
W_{2k} (\tau,  h),
\quad  h \in (0, 2\pi/n].
$$
}
\vskip .3cm
{\it Proof.} Denote by $\chi_h^r(x)$ the convolution power
of the normed characteristic function of the interval $[-h/2,
h/2]$.
$$
\chi_h^r(x):=(\chi_h*\chi_h^{r-1})(x), \quad  \int_T \chi_h(t) \, dt
=1.
$$
Note, that in these notations we have  $\phi_h(x)=\chi_h^2(x)$.
We can use the standard integral representation for the difference

$$ \widehat \Delta_t^r f(x) = t^r (D^r(f)*\chi_t^r)(x). $$
\vskip .3cm
\noindent
The Bernstein--Nikolsky--Stechkin inequality for  $ t \in (0, 2\pi/n)$
is equaivalent (see \cite{sbs}~) to the following inequality for $\tau \in T_n $,
$\| \tau \| = \tau(x_0) = 1$:
\vskip .3cm
\noindent
$$
(\tau*\chi_t^r)(x_0)\ge (c_n*\chi_t^r)(0).
$$
\vskip .3cm
\noindent
Therefore, after multiplication of last inequality (for $D^{2k} \tau$) by
$t^{2k}
\chi_h^2(t)$ and integration on  $t$ we get
$$
\int_T \widehat\Delta_t^{2k} \tau(x_0) \chi_h^2 (t) \, dt
=\int_T t^{2k} (D^{2k} \tau*\chi_t^{2k})(x_0) \chi_h^2 (t) \, dt
$$
$$
\ge \int_T
(c_n*\chi_t^{2k})(0)t^{2k} \chi_h^2 (t) \, dt, \quad  \| D^{2k} \tau \| = D^{2k} \tau (x_0)=1.
$$
\qed

\vskip .3cm

\newpage
\begin{center}
{\bf 4 Favard's operators and Jackson--Stechkin theorem}
\end{center}

\vskip .2cm

We shall call an operator  $A_{n,r}$   Favard's operator, if

$$
\| f - A_{n,r} (f) \| \le F_r n^{-r} \| D^r f  \|.
$$
Suppose that  $\tau_* \in T_{n-1}$ gives the estimate  \cite{fks, bk} :

$$
\| f - \tau_* \| \le c_{\alpha} W_{2k}(f, \alpha \pi/n),
\quad \alpha \in (1, n/k).
$$
\vskip .2cm
\noindent
Put $C_{k,\alpha}:= F_{2k} W_{2k}(c_n, \alpha \pi/n)^{-1},
\quad h_\alpha :=\alpha \pi/n$.
\vskip .2cm
\noindent
The Theorem 1 and  the Lemma  1 imply

$$
\| \tau_* - A_{n,2k} (\tau_*) \| \le F_{2k} n^{-2k} \| D^{2k} \tau_*  \| \le
C_{k,\alpha} W_{2k} (\tau_*, h_\alpha) \le
$$
\vskip .2cm
$$
C_{k,\alpha} \left( W_{2k}(f-\tau_*,h_\alpha) + W_{2k} (f, h_\alpha)
  \right) \le
 C_{k,\alpha} \left( 3 \| f - \tau_* \| +
W_{2k} (f, h_\alpha) \right).
$$
\vskip .2cm
\noindent
Thus

$$
\| f - A_{n,2k} (\tau_*) \| \le \| f - \tau_* \| + \| \tau_* - A_n (\tau_*) \|
\le
C(\alpha) W_{2k} (f, h_\alpha). \eqno (3)
$$
\vskip .2cm
\noindent
The methods of the paper  \cite{fks} allow us to obtain the following
estimates (see \cite{bk}):

$$
C_{k,\alpha} \le  \mathcal{K}_{2k} \,  c_\alpha, \quad  c_\alpha \le
\sec(\pi/(2\alpha)), \quad \mathcal{K}_{r}=\frac{4}{\pi} \sum_{j=-\infty}^\infty (4j+1)^{-r-1}.
$$
Therefore, we have
$$
C(\alpha) \le c_\alpha (1 + 3 C_{k, \alpha}) + C_{k,\alpha} \le c
\, (\alpha-1)^{-2}.
$$
\vskip .2cm
\noindent
In the case $\alpha = 1$ the estimates of the constants are the following (see \cite{bk} ):

$$
C_{k,1} \le  \mathcal{K}_{2k} \, c_1, \quad  c_1 \le c \, \sqrt {k},
$$
and
$$
C(1) = O(k).
$$

In the case of approximation of the periodic functions by
periodic smooth  splines the Favard's type estimates take place \cite{tih}
with

$$
F_r = \mathcal{K}_{r}.
$$

Thus, we have the Jackson--Stechkin inequality for approximation by periodic splines with
        best (respect to $r$) constants.

Rewrite the inequality (3) for approximation by Favard's operators in terms of
best approximations and standard moduli of smoothness. Let
$E_{n-1}^F (f)$ be the best uniform approximation of continuous
periodic function by Favard's  operators.

\newpage

{\bf Theorem  2.} For $r \in N$

$$
E_{n-1}^F(f)  \le c \, \max \left((\alpha-1)^{-2}, 1 \right) \, \sqrt{r} \  2^{-r}  \  \omega_r
(f, \alpha \pi/n), \quad \alpha \in ( 1, 2n/r). \eqno (4)
$$

$$
E_{n-1}^F (f)  \le c \  r^{3/2}  \  2^{-r}  \omega_r(f, \pi/n).
$$

\vskip .2cm

The sharpness of (4) with respect to order $r$ gives, for
example, the periodic step--functions: $ \mbox{\,sign \,}( c_n (t) )$.

\vskip .3cm
\small
\begin{tabular}[t]{l}
A.G. Babenko \\ Institute of Mathematics and Mechanics \\ Ural
Branch of the Russian Academy of Sciences \\ 16, S.Kovalevskoi
Str. \\ Ekaterinburg, 620219 \\ Russia \\
\textrm{babenko@imm.uran.ru}
\end{tabular}
 \hfill
\vskip .3cm
\begin{tabular}[t]{l}
Yuri Kryakin \\ Institute of Mathematics \\ University of Wroclaw
\\ Plac Grunwaldzki 2/4 \\ 50-384 Wroclaw \\ Poland \\
\textrm{kryakin@gmail.com}
\end{tabular}
\end{document}